\UseRawInputEncoding
\documentclass[11pt,twoside]{article}
\usepackage{amsmath, amsthm, amsfonts, amssymb}
\usepackage{mathrsfs}
\usepackage[colorlinks,linkcolor=blue,anchorcolor=blue,citecolor=blue,urlcolor=black]{hyperref}

\allowdisplaybreaks

\newfam\msbfam
\font\tenmsb=msbm5    \textfont\msbfam=\tenmsb \font\sevenmsb=msbm5
\scriptfont\msbfam=\sevenmsb \font\fivemsb=msbm5
\scriptscriptfont\msbfam=\fivemsb

\newfam\bigfam
\font\tenbig=msbm5 scaled \magstep2   \textfont\bigfam=\tenbig
\font\sevenbig=msbm7 scaled \magstep2 \scriptfont\bigfam=\sevenbig
\font\fivebig=msbm5 scaled \magstep2
\scriptscriptfont\bigfam=\fivebig

\def\dint{\displaystyle\int}

\def\dprod{\displaystyle\prod}
\def\dsup{\displaystyle\sup}
\def\dinf{\displaystyle\inf}

\makeatletter

\begin{document}
	
\title{\bf Multiple Weighted Estimates for Multilinear Commutators of Multilinear Singular Integrals with Generalized Kernels}
	
\author{\bf LIWEN GAO, YAN LIN AND SHUHUI YANG}
		
	\date{}
	\maketitle
	\renewcommand{\thefootnote}{}
	

\footnote{This work was financially supported by the National Natural Science Foundation of China (No. 12071052).}
	\begin{minipage}{13.5cm}
		
		 {\quad ABSTRACT. In this paper, the weighted $L^{p}$ boundedness of multilinear commutators and iterated commutators of multilinear singular integral operators with generalized kernels is established, where the weight is multiple weight. Our results are generalizations of the corresponding results for multilinear singular integral operators with standard kernels and Dini kernels.} 

	\end{minipage}
{\centering\section{\hspace{-0.3cm}{\bf} Introduction \label{s1}}}

  The multilinear singular integral operator theory is the most brilliant part of modern harmonic analysis. It originated from the work \cite{cm1}-\cite{cm3} of Meyer and Coifman in the 1970s. The theory mainly study the operator as follows:

\medskip
\noindent{{\bf Definition 1.1}}. Suoppse $m \in \mathbb{N}^{+}$ and a
function $K(y_{0},y_{1},\dots,y_{m})$ is defined apart from the diagonal $y_{0}=y_{1}=\dots=y_{m}$
in $(\mathbb{R}^{n})^{m+1}$. $T$ is an $m$-linear
operator defined on product of test functions that relates to $K$,
then its integral representation is as follow:
$$T(f_{1},\dots,f_{m})(x)=\dint_{\mathbb{R}^{n}}\cdots\dint_{\mathbb{R}^{n}}
K(x,y_{1},\dots,y_{m})\dprod_{j=1}^{m}f_{j}(y_{j})dy_{1}\cdots
dy_{m}, \eqno{(1.1)}$$ where $f_{j} (j=1,\dots,m)\in C_{0}^{\infty}(\mathbb{R}^{n})$ and $x\notin
\cap_{j=1}^{m}\textrm{supp} f_{j}$.

In particular, we call $K$ a standard $m$-linear Calder\'{o}n-Zygmund kernel if it has the following size estimates:
$$|K(y_{0}, y_{1}, \dots, y_{m})|\leq \frac{A}{(\sum_{k,\, l=0}^{m}
	|y_{k}-y_{l}|)^{mn}}, \eqno{(1.2)}$$ for some $A>0$ and all
$(y_{0}, y_{1}, \dots, y_{m})\in (\mathbb{R}^{n})^{m+1}$ are away from the diagonal,
and smooth estimates:
$$|K(y_{0}, \dots, y_{j}, \dots, y_{m})-K(y_{0}, \dots, y_{j}', \dots, y_{m})|
\leq \frac{C |y_{j}-y_{j}'|^{\varepsilon}}{(\sum_{k,\, l=0}^{m}
	|y_{k}-y_{l}|)^{mn+\varepsilon}},\eqno{(1.3)}$$ for some
$\varepsilon>0$, whenever $0\leq j\leq m$ and $|y_{j}-y_{j}'|\leq
\frac{1}{2}\max_{0\leq k\leq m} |y_{j}-y_{k}|$.

\medskip
\noindent{{\bf Definition 1.2}}. Let $T$ be an $m$-linear operator defined by (1.1) with a standard $m$-linear Calder\'{o}n-Zygmund kernel $K$,  we say that $T$ is a
 standard $m$-linear Calder\'{o}n-Zygmund operator if it satisfies either of the following condition.
 
    Given a group of numbers $t_{1},t_{2},\dots, t_{m}$
 meeting $1\leq t_{1},t_{2},\dots, t_{m}\leq \infty$ and
    $1/t=1/t_{1}+\dots +1/t_{m}$,
    
(1) $T$ maps $L^{t_{1},1}\times \dots \times L^{t_{m},1}$ into $L^{t,\infty}$ if $t>1$,

(2) $T$ maps $L^{t_{1},1}\times \dots \times L^{t_{m},1}$ into $L^{1}$ \quad if $t=1$,\\ 
where the signs $L^{t_{1},1},L^{t_{2},1},\dots,L^{t_{m},1}$ and $L^{t_{m},\infty}$ all stand for Lorentz spaces.

 The initial interest in the study of commutators is related to the generalization of the classical factorization theorem for Hardy space. Then P\'{e}rez and Torres generalized the definition to the commutator generated by multilinear singular integral operator and BMO function in \cite{pt}. 
 Now we recall the definition of the multilinear commutators and iterated commutators as follows. 

\medskip
\noindent{{\bf Definition 1.3}}.
Let $\vec{b}=(b_{1},\dots,b_{m})$ be a family of locally integrable functions. The $m$-linear commutator generated by $\vec{b}$ and the $m$-linear operator $T$ is defined by
$$T_{\vec{b}}(f_{1},\dots,f_{m})=\sum_{j=1}^{m}T_{\vec{b}}^{j}(\vec{f}),$$
whenever each term is the commutator of $b_{j}$  and the $j$th entry of $T$, that is
$$T_{\vec{b}}^{j}=b_{j}T(f_{1},\dots,f_{j},\dots,f_{m})
-T(f_{1},\dots,b_{j}f_{j},\dots,f_{m}).$$ 

\medskip
\noindent{{\bf Definition 1.4}}.
Suppose that $T$ is an $m$-linear singular integral operator and $\vec{b}=(b_{1},\dots,b_{m})$
is a family of locally integrable functions, then the $m$-linear iterated commutator generated by $T$ and $\vec{b}$ is defined to be
$$T_{\prod\vec{b}}(f_{1},\dots,f_{m})=[b_{1},[b_{2},\dots,[b_{m-1},[b_{m},T]_{m}]_{m-1} \dots ]_{2}]_{1} (\vec{f}),$$
where $\vec{f}=(f_{1},\dots,f_{m}).$
If $T$ is linked in the usual way to the kernel $K$ studied in this paper, then we can write
\begin{eqnarray*}
	T_{\prod\vec{b}}(f_{1},\dots,f_{m})(x)&=&\dint_{(\mathbb{R}^{n})^{m}}
	\prod_{j=1}^{m}\biggl(b_{j}(x)-b_{j}(y_{j}) \biggr) 
	K(x,y_{1},\dots,y_{m})f_{1}(y_{1})\cdots f_{m}(y_{m})\\
	&&\quad\times dy_{1} \cdots dy_{m}. 
\end{eqnarray*}

In 2002, Grafakos-Kalton in \cite{gk} and Grafakos-Torres in \cite{gt1}-\cite{gt3} gave a systematic statement of the theory. In \cite{gt1}, the authors showed the boundedness of classical singular integral operator. Grafakos-Martell made some work about weighted estimates for standard $m$-linear Calder\'{o}n-Zygmund operator in \cite{gm}. The weight that they discuss is the Muckenhoupt weight. 

They obtained that $T$ is bounded from
$L^{p_{1}}(w_{1})\times \dots \times L^{p_{m}}(w_{m})$ to $L^{p}(w)$ for any $1<p_{1},\dots,p_{m}<\infty$ and $1/p=1/p_{1}+\dots+1/p_{m}$.

In \cite{gt3}, Grafakos and Torres put forwards a question:" Is there a multiple weight theory? The most appropriate multilinear maximal function or multiple weights to work with in this direction are not yet clear." Therefore, to solve the question the multiple weight was given by Lerner, et al in \cite{loptt}. The definition of the multiple weight is as follows.

\medskip
\noindent{{\bf Definition 1.5}}.
Let $1\leq p_{1},\dots,p_{m}<\infty$. Given $\vec{w}=(w_{1},\dots,w_{m})$, set 
$$v_{\vec{w}}=\prod_{j=1}^{m}w_{j}^{\frac{p}{ p_{j} } }.$$
$\vec{w}$ satisfies the $A_{\vec{P}}$ condition if
$$\dsup_{Q}\biggl(\dfrac{1}{|Q|}\dint_{Q}v_{\vec{w}}\biggr)^{\frac{1}{p}}\prod_{j=1}^{m}\biggl(\dfrac{1}{|Q|}\dint_{Q} w{j}^{1-p_{j}'}\biggr)^{\frac{1}{p_{j}'}}<\infty.$$
When $p_{j}=1$, the condition $\biggl(\dinf_{Q}w_{j}\biggr)^{-1}$takes place of the  condition $\biggl(\dfrac{1}{|Q|}\dint_{Q} w{j}^{1-p_{j}'}\biggr)^{\frac{1}{p_{j}'}}$.

In \cite{loptt}, the authors not only put forwards the definition of the multiple weight, but also presented a new maximal function
$$\mathcal{M}(\vec{f})(x)=\dsup_{Q\ni x}\prod_{j=1}^{m}\dfrac{1}{|Q|}\dint_Q|f_{j}(y_{j})|dy_{j},$$
and denoted $\mathcal{M}_{q}(\vec{f})$ as
$$\mathcal{M}_{q}(\vec{f})(x)=\dsup_{Q\ni x}\biggl( \prod_{j=1}^{m}\dfrac{1}{|Q|}\dint_Q|f_{j}(y_{j})|^{q}dy_{j}
\biggr)^{\frac{1}{q}}.$$

This kind of maximal function will become an key tool in our paper in order to match the multiple weight. In \cite{loptt}, by using the key tools the authors obtained a range of results.
Set $\vec{w}$ under the assumption of Definition 1.5 and with regard to the standard multilinear Calder\'{o}n-Zygmund operator $T$, for any $1< p_{1},\cdots,p_{m}<\infty$ with $1/p=1/p_{1}+\cdots+1/p_{m}$, $T$ and $T_{\vec{b}}$ are both bounded from $L^{p_{1}}(w_{1})\times \cdots \times L^{p_{m}}(w_{m})$ into $L^{p}(v_{\vec{w}})$.

In 2014, Zhang-Lu introduced multilinear singular integral operators with kernels of Dini's type in \cite{lz2}. The operators that they studied weakened the condition (1.3) in the foundation of standard multilinear Calder\'{o}n-Zygmund operator into the following condition,
\begin{align*}
	&|K(y_{0}, \dots, y_{j}, \dots, y_{m})-K(y_{0}, \dots, y_{j}', \dots, y_{m})|
	\leq  \frac{C} {(|y_{0}-y_{1}|+\cdots+|y_{0}-y_{m}|)^{mn}}\\
	 &w\biggl(\frac{|y_{j}-y_{j}'|}{|y_{0}-y_{1}|+\cdots+|y_{0}-y_{m}|}\biggr) \tag{1.4},
\end{align*}
whenever $|y_{j}-y_{j}'|\leq\frac{1}{2}\max_{1\leq j\leq m} |y_{0}-y_{j}|$ and $w(t)$ is non-negative and non-decreasing function on $\mathbb{R}^{+}$. 

We easily know that the condition (1.4) can become the condition (1.3) by letting $w(t)=t^{\varepsilon}$. Thus the standard multilinear Calder\'{o}n-Zygmund operator is a special case of 
multilinear singular integral operators with kernels of Dini's type.

In \cite{lz2}, Zhang and Lu mainly studied the multilinear singular integrals with kernels of Dini's type and the multilinear commutators with multiple weight.

In \cite{lx}, Lin-Xiao gave a kind of multilinear singular integrals with generalized kernels which changed the condition (1.4) into a weaker condition.

For  any positive integers $k_{1},\dots,k_{m}$,
\begin{align*}
&\biggl(\int_{2^{k_{m}}|y_{0}-y_{0}'|\leq |y_{m}-y_{0}| < 2^{k_{m}+1}|y_{0}-y_{0}'|} \dots \int_{2^{k_{1}}|y_{0}-y_{0}'|\leq |y_{m}-y_{0}|< 2^{k_{1}+1}|y_{0}-y_{0}'|} \\
&|K(y_{0}, y_{1}, \dots, y_{m})-K(y_{0}', y_{1}, \dots, y_{m})|^{q} dy_{1}\dots dy_{m} \biggr)^{\frac{1}{q}} \\
&\leq C|y_{0}-y_{0}'|^{- \frac{mn}{q'}}\prod_{i=1}^{m}C_{k_{i}}
2^{-\frac{n}{q'}k_{i}}, \tag{1.5}
\end{align*}
where $1<q<\infty$, $C_{k_{i}}$ is a positive constant about $k_{i}$ and $q'$ is the conjugate exponent of $q$.

Motivated by \cite{gm,loptt,lx,lz1,lz2}, we now discuss the weighted boundedness of multilinear commutators and iterated commutators of multilinear singular integral operators with generalized kernels, where the weight is multiple weight. However, the condition (1.5) can not be suit for the multiple weight condition. This fact forces us to search for another kind of generalized kernels as follows.
For $k\in \mathbb{N}$,
\begin{align*}
	&\biggl(\int_{2^{k}\sqrt{m}|y_{0}-y_{0}'|\leq |\vec{y}-\vec{y_{0}}| <2^{k+1}\sqrt{m}|y_{0}-y_{0}'|}  |K(y_{0},y_{1}, \dots, y_{m})-K(y_{0}',y_{1}, \dots, y_{m})|^{q} dy_{1}\dots dy_{m}\biggr) ^{\frac{1}{q}} \\
	&\leq CC_{k}|y_{0}-y_{0}'|^{- \frac{mn}{q'}}2^{-\frac{kmn}{q'}}\tag{1.6},
\end{align*}
whenever $\vec{y_{0}}=(y_{0},y_{0},\dots,y_{0})$,
$\vec{y}=(y_{1},y_{2},\dots,y_{m})$, $(q,q')$ is a fixed pair of positive numbers satisfying $\frac{1}{q}+\frac{1}{q'}=1$ and $C_{k}$ is a positive constant about $k$.

\medskip
\noindent{{\bf Remark 1.6}}.
By calculation, we can find that condition (1.4) implies condition (1.6) by putting $C_{k}=w(2^{-k})$. According to the contents that we mentioned in the foregoing, the work in \cite{lz1} includes the classical case. So we can deduce naturally that our work contains both the cases of standard kernel and Dini type kernel.

  This paper is organzed as follows. In Section 2, we list the necessary definition and lemmas. In Secion 3, we establish the weighted estimates of $T_{\vec{b}}$ and $T_{\prod\vec{b}}$ by their sharp estimation. In Section 4, we will give the proof of all theorems.\par
 \centerline{} 
{\centering\section{\hspace{-0.3cm}{\bf} Necessary Definition and Lemma \label{s2}}}
Meanwhile, we also give some definition and lemmas used in this paper.

\medskip
\noindent{{\bf Definition 2.1}}.
Fix positive integers $k$ and $m$ satisfying $1\leq k<m$ and we suppose $C_{k}^{m}$ be a family of all finite subsets $\sigma=\{\sigma(1),\dots ,\sigma(k)\}$ of $\{1,2,\dots,m\}$ with $k$ different elements. If $j<l$, then $\sigma(j)<\sigma(l)$. For any 
$\sigma\in C_{k}^{m}$, let $\sigma'=\{1,2,\dots,m\}\backslash\sigma$ be the complementary sequence. In particular, $C_{0}^{m}=\varnothing$. For an $m$-tuple $\vec{b}$ and $\sigma\in C_{k}^{m}$, the $k$-tuple
$\vec{b_{\sigma}}=(b_{\sigma(1)},\dots,b_{\sigma(k)})$ is a finite subset of $\vec{b}=(b_{1},\dots,b_{m})$.
Given $T$ the $m$-linear operator and $\sigma\in C_{j}^{m}$, the multilinear iterated
commutator is given by
$$T_{\prod\vec{b_{\sigma}}}(f_{1},\dots,f_{m})=[b_{\sigma(1)},[b_{\sigma(2)},\dots,[b_{\sigma(j-1)},[b_{\sigma(j)},T]_{\sigma(j)}]_{\sigma(j-1)} \dots ]_{\sigma(2)}]_{\sigma(1)} (\vec{f}).$$
According to (1.1), we can also represent the multilinear iterated commutator by the form.
\begin{eqnarray*}
	T_{\prod\vec{b_{\sigma}}}(f_{1},\dots,f_{m})(x)&=&\dint_{(\mathbb{R}^{n})^{m}}
	\prod_{j=1}^{m}\biggl(b_{\sigma(i)}(x)-b_{\sigma(i)}(y_{\sigma(i)}) \biggr) 
	K(x,y_{1},\dots,y_{m})\\
	&&\times f_{1}(y_{1})\cdots f_{m}(y_{m})dy_{1} \cdots dy_{m} .
\end{eqnarray*}
Obviously, $T_{\prod\vec{b_{\sigma}}}=T_{\prod\vec{b}}$ when 
$\sigma=\{1,2,\dots,m\}$ and $T_{\prod\vec{b_{\sigma}}}=T_{b_{j}}^{j}$
when $\sigma=\{j\}.$

 \medskip
 \noindent{{\bf Definition 2.2}}.
 The Hardy-Littlewood maximal operator $M$ is defined by
 $$Mf(x)=\dsup_{B\ni x}\dfrac{1}{|B|}\dint_B|f(y)|dy,$$
 where the supremum is taken over all balls which contain $x$. We can also define the operator $M_{s}(f)(x)$ by $M_{s}(f)=[M(|f|^{s})]^{1/s}$, $s>0.$
 The sharp maximal operator $M^{\sharp}$ is defined by
$$M^{\sharp}(f)(x)=\dsup_{B\ni
	x}\dfrac{1}{|B|}\dint_{B}|f(y)-f_{B}|dy\sim \dsup_{B\ni x} \dinf_{a\in
	\mathbb{C}}\dfrac{1}{|B|}\dint_{B}|f(y)-a|dx.$$
We denote the $l$-sharp maximal operator by
$M^{\sharp}_{l}(f)=[M^{\sharp}(|f|^{l})]^{1/l}$, $l>0.$

\medskip
\noindent{{\bf Lemma 2.3}}(see \cite{loptt}). {\it Suppose that $p$, $q$ satisfy $0<p<q<\infty$. Then there is a positive constant $C=C_{p,q}$ so as to for any measurable function $f$, the following equality
	$$|Q|^{-1/p}\|f\|_{L^{p}(Q)}\leq C|Q|^{-1/q}\|f\|_{L^{q,\infty}(Q)}$$}
holds for any $Q$.

\medskip
\noindent{{\bf Lemma 2.4}}(see \cite{ll}). {\it
	Suppose that $f$ is a function in $BMO(\mathbb{R}^{n})$, $1\leq p<\infty$,
	$x\in \mathbb{R}^{n}$, and $r_{1}, r_{2}>0$. Then
	there is a positive constant $C$ independent of $f$, $x$, $r_{1}$ and $r_{2}$ such that 
	$$\biggl(\frac{1}{|B(x, r_{1})|}\dint_{B(x,\, r_{1})}|f(y)-f_{B(x,\,
		r_{2})}|^{p} dy\biggr)^{1/p}\leq
	C\biggl(1+\bigl|\ln\frac{r_{1}}{r_{2}}\bigr|\biggr)\|f\|_{
		BMO}.$$ }
	
\medskip
\noindent{{\bf Lemma 2.5}}(see \cite{loptt}, \cite{fs}). {\it
	Let $0< p,\delta<\infty$ and $w\in A_{\infty}$. Then there exists a constant $C>0$ depending
	only on the $A_{\infty}$ constant of $w$ such that
	$$\dint_{\mathbb{R}^{n}}[M_{\delta}(f)(x)]^{p}w(x)dx\leq C\dint_{\mathbb{R}^{n}}[M_{\delta}^{\sharp}(f)(x)]^{p}w(x)dx,$$
	for every function $f$ such that the left-hand side is finite.}

\medskip
\noindent{{\bf Lemma 2.6}}(see \cite{loptt}). {\it
Let $\vec{w}=(w_{1},\dots,w_{m})$ and $1\leq p_{1},\cdots,p_{m}<\infty$. Then $w\in A_{\vec{P}}$ if and only if
\begin{eqnarray*}
\left\{
\begin{array}{lll}
w_{j}^{1-p_{j}'} \in A_{mp_{j}'}, &j=1,2,\dots,m,\\
v_{\vec{w}} \in A_{mp}.
\end{array}
\right.
	\end{eqnarray*}
The condition	$w_{j}^{1-p_{j}'} \in A_{mp_{j}'}$ in the case $p_{j}=1$ is replaced by $w_{j}^{\frac{1}{m}}\in A_{1}$.}

\medskip
\noindent{{\bf Lemma 2.7}}(see \cite{loptt}). {\it
Let $1<p_{j}<\infty$, $j=1,2,\dots ,m$, and $1/p=1/p_{1}+\cdots+1/p_{m}$. Then the inequality	
$$\| \mathcal{M}(\vec{f})\|_{ L^{p}(v_{\vec{w}}) }
\leq C \prod_{j=1}^{m}\|f_{j}\|_{L^{p_{j}}(w_{j})}$$
holds for any $\vec{f}$ if and only if $\vec{w}$ is content with  $A_{\vec{P}}$ condition. 
	
}

\medskip
\noindent{{\bf Lemma 2.8}}(see \cite{yll}). {\it Let $m\geq 2$ and $T$ be an
	$m$-linear singular integral operator defined by (1.1) with generalized kernel satisfying (1.6) and $\sum_{k=1}^{\infty}C_{k}<\infty$. Suppose for fixed $1\leq r_{1},\cdots,r_{m}\leq q'$ with $1/r=1/r_{1}+\cdots+1/r_{m}$,
	 $T$ is bounded from $L^{r_{1}}\times \cdots \times L^{r_{m}}$ into $L^{r,\infty}$. If $0<\delta<1/m$, then we have
	$$M^{\sharp}_{\delta}(T(\vec{f}))(x)\leq C\mathcal{M}_{q'}(\vec{f})(x),$$
	for all $m$-tuples $\vec{f}=(f_{1},\dots,f_{m})$ of bounded measurable functions with compact support. }\par
\centerline{} 
{\centering\section{\hspace{-0.3cm}{\bf} Main results \label{s4}}}

\medskip
\noindent{{\bf Theorem 3.1.}} {\it Let $m\geq 2$ and $T$ be an
	$m$-linear singular integral operator defined by (1.1) with generalized kernel satisfying (1.6) and $\sum_{k=1}^{\infty}kC_{k}<\infty$. Suppose for fixed $1\leq r_{1},\cdots,r_{m}\leq q'$ with $1/r=1/r_{1}+\cdots+1/r_{m}$, $T$ is bounded from $L^{r_{1}}\times \cdots \times L^{r_{m}}$ into $L^{r,\infty}$. If $\vec{b}\in BMO^{m} $, $0<\delta<1/m$, $\delta<\varepsilon_{0}<\infty$ and $q'<s<\infty$,
	then for all m-tuples $\vec{f}=(f_{1},\cdots,f_{m})$ of bounded measurable functions with compact support, there is 
	$$M^{\sharp}_{\delta}(T_{\vec{b}}(\vec{f}))(x)\leq C\| \vec{b} \|_{BMO^{m}} \biggl(M_{\varepsilon_{0}}(T(\vec{f}))(x)+\mathcal{M}_{s}(\vec{f})(x)\biggr),$$
	where $\|\vec{b} \|_{BMO^{m}}=\max_{1\leq j\leq m}\|b_{j} \|_{BMO}$.}

	\medskip
	\noindent{{\bf Theorem 3.2.}} {\it Let $m\geq 2$ and $T$ be an $m$-linear singular integral operator defined by (1.1) with generalized kernel satisfying (1.6) and $\sum_{k=1}^{\infty}kC_{k}<\infty$. Suppose for fixed $1\leq r_{1},\cdots,r_{m}\leq q'$ with $1/r=1/r_{1}+\cdots+1/r_{m}$,
		$T$ is bounded from $L^{r_{1}}\times \cdots \times L^{r_{m}}$ into $L^{r,\infty}$. If $\vec{b}\in BMO^{m} $, then for $q'< p_{1},\cdots,p_{m}<\infty$ with $1/p=1/p_{1}+\cdots+1/p_{m}$, $T_{\vec{b}}$ is bounded from $L^{p_{1}}(w_{1})\times \cdots \times L^{p_{m}}(w_{m})$ into $L^{p}(v_{\vec{w}})$, where  $\vec{w}=(w_{1},\dots,w_{m})\in A_{\vec{P}/q'}$ and $v_{\vec{w}}=\prod_{j=1}^{m}w_{j}^{\frac{p}{ p_{j} } }$.}

\medskip
\noindent{{\bf Theorem 3.3.}} {\it Let $m\geq 2$ and $T$ be an
	$m$-linear singular integral operator defined by (1.1) with generalized kernel satisfying (1.6) and $\sum_{k=1}^{\infty}k^{m}C_{k}<\infty$. Suppose for fixed $1\leq r_{1},\cdots,r_{m}\leq q'$ with $1/r=1/r_{1}+\cdots+1/r_{m}$,
	$T$ is bounded from $L^{r_{1}}\times \cdots \times L^{r_{m}}$ into $L^{r,\infty}$. If $\vec{b}\in BMO^{m} $, $0<\delta<1/m$, $\delta<\varepsilon_{0}<\infty$ and $q'<s<\infty$,
	then for all m-tuples $\vec{f}=(f_{1},\cdots,f_{m})$ of bounded measurable functions with compact support,}
\begin{eqnarray*}
M^{\sharp}_{\delta}(T_{\prod\vec{b}}(\vec{f}))(x)&\leq& C\prod_{j=1}^{m}\| b_{j} \|_{BMO} \biggl(M_{\varepsilon_{0}}(T(\vec{f}))(x)+\mathcal{M}_{s}(\vec{f})(x)\biggr)\\
&&+C\sum_{j=1}^{m-1}\sum_{\sigma\in C_{j}^{m}}\prod_{i=1}^{j}\| b_{\sigma(i)} \|_{BMO} M_{\varepsilon_{0}}(T_{\prod\vec{b}_{\sigma' }}(\vec{f}))(x).
\end{eqnarray*}
\medskip
\noindent{{\bf Theorem 3.4.}} {\it  Let $m\geq 2$ and $T$ be an
	$m$-linear singular integral operator defined by (1.1) with generalized kernel satisfying (1.6) and $\sum_{k=1}^{\infty}k^{m}C_{k}<\infty$. Suppose for fixed $1\leq r_{1},\cdots,r_{m}\leq q'$ with $1/r=1/r_{1}+\cdots+1/r_{m}$,
	$T$ is bounded from $L^{r_{1}}\times \cdots \times L^{r_{m}}$ into $L^{r,\infty}$. If $\vec{b}\in BMO^{m} $, $q'<p_{1},\cdots,p_{m}<\infty$ with $1/p=1/p_{1}+\cdots+1/p_{m}$, then $T_{\prod \vec{b}}$ is bounded from $L^{p_{1}}(w_{1})\times \cdots \times L^{p_{m}}(w_{m})$ into  $L^{p}(v_{\vec{w}})$, where $\vec{w}=(w_{1},\dots,w_{m})\in A_{\vec{P}/q'}$ and $v_{\vec{w}}=\prod_{j=1}^{m}w_{j}^{\frac{p}{ p_{j} } }$.}\par
\centerline{} 
{\centering\section{\hspace{-0.3cm}{\bf} Proof of the main results \label{s5}}}

\medskip
\noindent{{\it Proof of Theorem $3.1$.}}
Fix $x\in\mathbb{R}^{n}$. For any ball $B\ni x$, we first consider the commutator.
 $$T_{b_{1}}^{1}(f_{1},\dots,f_{m})(z)=b_{1}(z)T(f_{1},\dots,f_{m})(z)-T(b_{1}f_{1},\dots,f_{m})(z).$$
Let $B^{*}=16\sqrt{m}B$, then for any $z\in B$,
 $$T_{b_{1}}^{1}(f_{1},\dots,f_{m})(z)=(b_{1}(z)-b_{1B^{*}})T(f_{1},\dots,f_{m})(z)-T((b_{1}-b_{1B^{*}})f_{1},\dots,f_{m})(z),$$
 where $b_{1B^{*}}=\dfrac{1}{|B^{*}|}\dint_{B^{*}}b_{1}(z)dz.$
 Since $0<\delta<1$, then for any $c\in \mathbb{C}$, we have
 \begin{eqnarray*}
 	&&\biggl(\frac{1}{|B|}\dint_{B}||T_{b_{1}}^{1}(\vec{f})(z)|^{\delta}-|c|^{\delta}|dz\biggr)^{\frac{1}{\delta}}\\
 	&\leq& C\biggl(\frac{1}{|B|}\dint_{B}|(b_{1}(z)-b_{1B^{*}})T(\vec{f})(z)|^{\delta}dz\biggr)^{\frac{1}{\delta}}\\
 	&&+C\biggl(\frac{1}{|B|}\dint_{B}|T((b_{1}-b_{1B^{*}})f_{1},\dots,f_{m})(z)-c|^{\delta}|dz\biggr)^{\frac{1}{\delta}}\\
 	&:=&I+II.
  \end{eqnarray*}

 We can find an $l$ such that $1<l<\min\{\frac{\varepsilon_{0}}{\delta}$, $\frac{1}{1-\delta}\}$, then $l\delta<\varepsilon_{0}$ and $l'\delta>1$. By H\"{o}lder's inequality and Lemma 2.4, we obtain
\begin{eqnarray*}
I&\leq& C\biggl(\frac{1}{|B|}\dint_{B}|b_{1}(z)-b_{1B^{*}}|^{\delta l'}dz\biggr)^{\frac{1}{\delta l'}}\biggl(\frac{1}{|B|}\dint_{B}|T(\vec{f})(z)|^{\delta l}dz\biggr)^{\frac{1}{\delta l}}\\
&\leq&C (1+\ln16\sqrt{m})\|b_{1}\|_{BMO}M_{\delta l}(T(\vec{f}))(x)\\
&\leq&C\|b_{1}\|_{BMO}M_{\varepsilon_{0}}(T(\vec{f}))(x).
 \end{eqnarray*}

 For each $j$, 
 $f_{j}=f_{j}^{0}+f_{j}^{\infty}$ where  
 $f_{j}^{0}=f_{j}\chi_{B^{*}}$. Then
$$
\prod_{j=1}^{m}f_{j}(y_{j})
=\prod_{j=1}^{m}f_{j}^{0}(y_{j})+\sum_{(\alpha_{1},\cdots,\alpha_{m})\in\Gamma}
f_{1}^{\alpha_{1}}(y_{1})\cdots f_{m}^{\alpha_{m}}(y_{m}),
$$
where $\Gamma=\{(\alpha_{1},\alpha_{2},\cdots,\alpha_{m}): \enspace\text{there is at least one}\enspace \alpha_{j}\neq0\}.$ Denote $\vec{z_{0}}=(z_{0},z_{0},\cdots,z_{0})$,
$\vec{y}=(y_{1},y_{2},\cdots,y_{m})$, the center of $B$ by $x_{0}$
and the radius of $B$ by $r_{B}$. We choose $z_{0}\in 4B\backslash3B$.
Let $$c=\sum_{(\alpha_{1},\cdots,\alpha_{m})\in\Gamma}T((b_{1}-b_{1B^{*}})f_{1}^{\alpha_{1}},f_{2}^{\alpha_{2}},\dots,f_{m}^{\alpha_{m}})(z_{0}).$$
Then 
\begin{eqnarray*}
	II&\leq&C\biggl( \frac{1}{|B|}\dint_{B}|T((b_{1}-b_{1B^{*}})f_{1}^{0},f_{2}^{0},\cdots,f_{m}^{0})(z)|^{\delta}dz \biggr)^{\frac{1}{\delta}}	\\
	&&+C\sum_{(\alpha_{1},\dots,\alpha_{m})\in\Gamma}\biggl( \frac{1}{|B|}\dint_{B}|T((b_{1}-b_{1B^{*}})f_{1}^{\alpha_{1}},f_{2}^
	{\alpha_{2}},\cdots,f_{m}^{\alpha_{m}})(z)\\
	&&-T((b_{1}-b_{1B^{*}})f_{1}^{\alpha_{1}},f_{2}^
	{\alpha_{2}},\cdots,f_{m}^{\alpha_{m}})(z_{0})|^{\delta}dz \biggr)^{\frac{1}{\delta}}\\
	&:=&II_{0}+\sum_{(\alpha_{1},\dots,\alpha_{m})\in\Gamma}II_{\alpha_{1},\cdots,\alpha_{m}}.
	\end{eqnarray*}

 Let $t=s/q'$, then it follows from $s>q'$ that $t>1$. Due to for any $j$, $r_{j}\leq q'$, we can get $r_{j}t\leq s.$
 By Lemma 2.3 and Lemma 2.4,
 we have
 \begin{eqnarray*}
  II_{0}&\leq& C|B|^{-1/\delta}\|T((b_{1}-b_{1B^{*}})f_{1}^{0},f_{2}^{0},\cdots,f_{m}^{0})\|_{L^{\delta}(B)} \\
    &\leq& C|B|^{-1/r}\|T((b_{1}-b_{1B^{*}})f_{1}^{0},f_{2}^{0},\cdots,f_{m}^{0})\|_{L^{r,\infty}(\mathbb{R}^{n})} \\
 	&\leq& C\biggl(\frac{1}{|B^{*}|}\dint_{B^{*}}  |(b_{1}(z)-b_{1B^{*}})f_{1}(y_{1})|^{r_{1}}dy_{1}\biggr)^{\frac{1}{r_{1} }}\prod_{j=2}^{m} \biggl(\frac{1}{|B^{*}|}\dint_{B^{*}}  |f_{j}(y_{j})|^{r_{j}}dy_{j}\biggr)^{\frac{1}{r_{j} }}\\
 	&\leq& C\biggl(\frac{1}{|B^{*}|}\dint_{B^{*}}  |b_{1}(z)-b_{1B^{*}}|^{r_{1}t'}dy_{1}\biggr)^{\frac{1}{r_{1} t'}}\biggl(\frac{1}{|B^{*}|}\dint_{B^{*}}  |f_{1}(y_{1})|^{r_{1}t}dy_{1}\biggr)^{\frac{1}{r_{1}t }}\\
 	&&\times \prod_{j=2}^{m} \biggl(\frac{1}{|B^{*}|}\dint_{B^{*}}  |f_{j}(y_{j})|^{r_{j}}dy_{j}\biggr)^{\frac{1}{r_{j}}}\\
 	&\leq& C\|b_{1}\|_{BMO}\prod_{j=1}^{m} \biggl(\frac{1}{|B^{*}|}\dint_{B^{*}} |f_{j}(y_{j})|^{s}dy_{j}\biggr)^{\frac{1}{s}}\\
 	&\leq& C \|b_{1}\|_{BMO}\mathcal{M}_{s}(\vec{f})(x).
 \end{eqnarray*}

For any $(\alpha_{1},\alpha_{2},\dots,\alpha_{m}) \in \Gamma $, we can find a $j\in\{1,\dots,m\}$ that makes $\alpha_{j}=\infty$. Then for any $ \vec{y}=(y_{1},y_{2},\dots,y_{m})\in {\rm supp}f_{1}^{\alpha_{1}}\times \dots\times {\rm supp}f_{m}^{\alpha_{m}}$ and $z\in B$, $|\vec{y}-\vec{z_{0}}|\geq2\sqrt{m}|z-z_{0}|$, and $2r_{B}\leq|z-z_{0}|\leq5r_{B}.$
Thus
 \begin{eqnarray*}
	II_{\alpha_{1},\alpha_{2},\dots,\alpha_{m}}
	&\leq&C\frac{1}{|B|}\dint_{B} \dint_{|\vec{y}-\vec{z_{0}}|\geq 2\sqrt{m}|z-z_{0}|} |K(z,\vec{y})-K(z_{0},\vec{y})||b_{1}(y_{1})-b_{1B^{*}}|\\
	&&\times\prod_{j=1}^{m}|f_{j}(y_{j})|d\vec{y}dz
	\\
	&\leq&C \frac{1}{|B|}\dint_{B}\sum_{k=1}^{\infty}\dint_{2^{k}\sqrt{m}|z-z_{0}|\leq|\vec{y}-\vec{z_{0}}|< 2^{k+1}\sqrt{m}|z-z_{0}|} |K(z,\vec{y})-K(z_{0},\vec{y})|\\
	&&\times|b_{1}(y_{1})-b_{1B^{*}}|\prod_{j=1}^{m}|f_{j}(y_{j})|d\vec{y}dz\\
	&\leq&C\frac{1}{|B|}\dint_{B} \sum_{k=1}^{\infty}\biggl(\dint_{2^{k}\sqrt{m}|z-z_{0}|\leq|\vec{y}-\vec{z_{0}}|< 2^{k+1}\sqrt{m}|z-z_{0}|} |K(z,\vec{y})-K(z_{0},\vec{y})|^{q}d\vec{y}\biggr)^{\frac{1}{q}}\\
	&&\times\biggl(\dint_{|\vec{y}-\vec{z_{0}}|< 2^{k+1}\sqrt{m}|z-z_{0}|}|b_{1}(y_{1})-b_{1B^{*}}|^{q'}\prod_{j=1}^{m}|f_{j}(y_{j})|^{q'}d\vec{y}\biggr)^{\frac{1}{q'}}dz\\
	&\leq&C\frac{1}{|B|}\dint_{B}
	\sum_{k=1}^{\infty}\biggl(\dint_{2^{k}\sqrt{m}|z-z_{0}|\leq|\vec{y}-\vec{z_{0}}|< 2^{k+1}\sqrt{m}|z-z_{0}|} |K(z,\vec{y})-K(z_{0},\vec{y})|^{q}d\vec{y}\biggr)^{\frac{1}{q}}\\
	&&\times\biggl(\frac{1}{|B(z_{0},2^{k+1}\sqrt{m}|z-z_{0}|)|}\dint_{|y_{1}-z_{0}|< 2^{k+1}\sqrt{m}|z-z_{0}|}
	|(b_{1}(y_{1})-b_{1B^{*}})\\
	&&\times f_{1}(y_{1})|^{q'}dy_{1}\biggr)^{\frac{1}{q'}}\\
	&&\times\prod_{j=2}^{m}\biggl(\frac{1}{|B(z_{0},2^{k+1}\sqrt{m}|z-z_{0}|)|}\dint_{|y_{j}-z_{0}|< 2^{k+1}\sqrt{m}|z-z_{0}|}
	|f_{j}(y_{j})|^{q'}dy_{j}\biggr)^{\frac{1}{q'}}\\
	&&\times(2^{k+1}\sqrt{m}|z-z_{0}|)^{ \frac{mn}{q'}}dz\\ 	
	&\leq&C\frac{1}{|B|}\dint_{B}
	\sum_{k=1}^{\infty}\biggl(\dint_{2^{k}\sqrt{m}|z-z_{0}|\leq|\vec{y}-\vec{z_{0}}|< 2^{k+1}\sqrt{m}|z-z_{0}|} |K(z,\vec{y})-K(z_{0},\vec{y})|^{q}d\vec{y}\biggr)^{\frac{1}{q}}\\
	&&\times\biggl(\frac{1}{|B(x_{0},2^{k+2}\sqrt{m}|z-z_{0}|)|}\dint_{|y_{1}-x_{0}|< 2^{k+2}\sqrt{m}|z-z_{0}|}
	|b_{1}(y_{1})-b_{1B^{*}}|^{q't'}dy_{1}\biggr)^{\frac{1}{q't'}}\\
	&&\times\biggl(\frac{1}{|B(z_{0},2^{k+1}\sqrt{m}|z-z_{0}|)|}\dint_{|y_{1}-z_{0}|< 2^{k+1}\sqrt{m}|z-z_{0}|}
	|f_{1}(y_{1})|^{q't}dy_{1}\biggr)^{\frac{1}{q't}}\\
	&&\times\prod_{j=2}^{m}\biggl(\frac{1}{|B(z_{0},2^{k+1}\sqrt{m}|z-z_{0}|)|}\dint_{|y_{j}-z_{0}|< 2^{k+1}\sqrt{m}|z-z_{0}|}
	|f_{j}(y_{j})|^{q't}dy_{j}\biggr)^{\frac{1}{q't}}\\
	&&\times(2^{k+1}\sqrt{m}|z-z_{0}|)^{ \frac{mn}{q'}}dz\\	
	&\leq&C\|b_{1}\|_{BMO}\frac{1}{|B|}\dint_{B}
	\sum_{k=1}^{\infty}|z-z_{0}|^{- \frac{mn}{q'}}C_{k}
	2^{-\frac{kmn}{q'}}k(2^{k+1}\sqrt{m}|z-z_{0}|)^{ \frac{mn}{q'}}\\
	&&\times\prod_{j=1}^{m}\biggl(\frac{1}{|B(z_{0},2^{k+1}\sqrt{m}|z-z_{0}|)|}\dint_{|y_{j}-z_{0}|< 2^{k+1}\sqrt{m}|z-z_{0}|}
	|f_{j}(y_{j})|^{s}dy_{j}\biggr)^{\frac{1}{s}}dz\\
	&\leq& C\|b_{1}\|_{BMO}\mathcal{M}_{s}(\vec{f})(x).
\end{eqnarray*}
Thus
$$
\sum_{(\alpha_{1},\dots,\alpha_{m})\in\Gamma}II_{\alpha_{1},\dots,\alpha_{m}}
\leq C\|b_{1}\|_{BMO}\mathcal{M}_{s}(\vec{f})(x).
$$
Through the above estimates we get the result 
\begin{eqnarray*}
M^{\sharp}_{\delta}(T_{\vec{b}}^{1}(\vec{f}))(x)&\leq& C\| b_{1} \|_{BMO} \biggl(M_{\varepsilon_{0}}(T(\vec{f}))(x)+\mathcal{M}_{s}(\vec{f})(x)\biggr)\\
&\leq& C\| \vec{b}\|_{BMO} \biggl(M_{\varepsilon_{0}}(T(\vec{f}))(x)+\mathcal{M}_{s}(\vec{f})(x)\biggr).
\end{eqnarray*}
Similarly, we can get to for any $j=1,2,\dots,m$,
$$
M^{\sharp}_{\delta}(T_{\vec{b}}^{j}(\vec{f}))(x)
\leq C\| \vec{b} \|_{BMO} \biggl(M_{\varepsilon_{0}}(T(\vec{f}))(x)+\mathcal{M}_{s}(\vec{f})(x)\biggr).
$$
So we have
\begin{eqnarray*}
	M^{\sharp}_{\delta}(T_{\vec{b}}(\vec{f}))(x)&=&M^{\sharp}_{\delta} \left(\sum_{j=1}^{m}T_{\vec{b}}^{j}(\vec{f})\right)(x)\\
	&\leq&
	\sum_{j=1}^{m}M^{\sharp}_{\delta}(T_{\vec{b}}^{j}(\vec{f}))(x)\\
	&\leq& C\| \vec{b} \|_{BMO} \biggl(M_{\varepsilon_{0}}(T(\vec{f}))(x)+\mathcal{M}_{s}(\vec{f})(x)\biggr),
\end{eqnarray*}
which completes the proof of the theorem.$\hfill\square$

\medskip
\noindent{{\it Proof of Theorem $3.2$.}}
Let $\vec{q}=\frac{\vec{P}}{q'}=(\frac{p_{1}}{q'},\dots,\frac{p_{m}}{q'})$. Due to $w\in A_{\vec{q}}$, by Lemma 2.6, each $\psi_{j}=w_{j}^{-\frac{1}{q_{j}-1}}$ belongs to $A_{\infty}$
where $q_{j}=\frac{p_{j}}{q'}$, $j=1,\dots,m$. By inverse H\"{o}lder's inequality, we can find constants $c_{j},t_{j}>1$ depending on the $A_{\infty}$ constant of $\psi_{j}$
such that
$$\biggl(\frac{1}{|B|}\int_{B}w_{j}^{-\frac{t_{j}}{q_{j}-1}}\biggr)^{\frac{1}{t_{j}}}\leq \frac{c_{j}}{|B|}\int_{B}w_{j}^{-\frac{1}{q_{j}-1}},$$
for any ball $B$. Pick a $d_{j}>1$ that makes 
$$\frac{t_{j}}{q_{j}-1}=\frac{1}{\frac{q_{j}}{d_{j}}-1}.$$
Then $q_{j}>d_{j}>1$, $j=1,\dots,m$.

Let $d=\min\{d_{1},\dots,d_{m}\}$ and $c=\max\{c_{1},\dots,c_{m}\}$. We have for $q_{0}=\frac{p}{q'}$
\begin{eqnarray*}
	&&\biggl(\frac{1}{|B|}\int_{B}v_{\vec{w}}\biggr)^{1/\frac{q_{0}}{d}}
	\prod_{j=1}^{m} \biggl(\frac{1}{|B|}\int_{B}w_{j}^{-\frac{1}{\frac{q_{j}}{d}-1}}\biggr)^{1-1/\frac{q_{j}}{d}}\\
	&=&\biggl(\frac{1}{|B|}\int_{B}v_{\vec{w}}\biggr)^{\frac{d}{q_{0}}}
	\prod_{j=1}^{m} \biggl(\frac{1}{|B|}\int_{B}w_{j}^{-\frac{1}{\frac{q_{j}}{d}-1}}\biggr)^{(\frac{q_{j}}{d}-1)\frac{d}{q_{j}}}\\
	&\leq&\biggl(\frac{1}{|B|}\int_{B}v_{\vec{w}}\biggr)^{\frac{d}{q_{0}}}\prod_{j=1}^{m} \biggl(\frac{1}{|B|}\int_{B}w_{j}^{-\frac{1}{\frac{q_{j}}{d_{j}}-1}}\biggr)^{(\frac{q_{j}}{d_{j}}-1)\frac{d}{q_{j}}}\\
	&=&\biggl(\frac{1}{|B|}\int_{B}v_{\vec{w}}\biggr)^{\frac{d}{q_{0}}}
	\prod_{j=1}^{m} \biggl(\frac{1}{|B|}\int_{B}w_{j}^{-\frac{t_{j}}{q_{j}-1}}\biggr)^{(\frac{q_{j}-1}{t_{j}})\frac{d}{q_{j}}}\\	
	&\leq&c^{dm}\biggl(\frac{1}{|B|}\int_{B}v_{\vec{w}}\biggr)^{\frac{d}{q_{0}}}
	\prod_{j=1}^{m} \biggl(\frac{1}{|B|}\int_{B}w_{j}^{-\frac{1}{q_{j}-1}}\biggr)^{(q_{j}-1)\frac{d}{q_{j}}}\\
	&\leq&c^{dm}[w]_{A_{\vec{q}}}^{d}.
\end{eqnarray*}

Thus $\vec{w}\in A_{\vec{q}/d}$. Let $s=q'd$, then we can find that $s>q'$ and $s<p_{j}$, $j=1,\dots,m$ from $d>1$. We can get the result
$$\vec{w}\in A_{\vec{q}/d}=A_{\vec{P}/s}.$$

Take $\varepsilon_{0}$, $\delta$ so that $0<\delta<\varepsilon_{0}<\frac{1}{m}.$
 We can deduce the following conclusion from Lemma 2.5, Theorem 3.1, Lemma 2.8, Lemma 2.7
and $\vec{w}\in A_{\vec{P}/s}$.
\begin{eqnarray*}
	\|T_{\vec{b}}(\vec{f})\|_{L^{p}(v_{\vec{w}})}
	&\leq&\|M_{\delta}(T_{\vec{b}}(\vec{f}))\|_{L^{p}(v_{\vec{w}})}\\
	&\leq&\|M_{\delta}^{\sharp}(T_{\vec{b}}(\vec{f}))\|_{L^{p}(v_{\vec{w}})}\\
	&\leq& C\| \vec{b} \|_{BMO^{m}} \|M_{\varepsilon_{0}}(T(\vec{f}))+\mathcal{M}_{s}(\vec{f})\|
	_{L^{p}(v_{\vec{w}})}\\
	&\leq& C\| \vec{b} \|_{BMO^{m}} (\|M_{\varepsilon_{0}}^{\sharp}(T(\vec{f}))\|_{L^{p}(v_{\vec{w}})}+\|\mathcal{M}_{s}(\vec{f})\|
	_{L^{p}(v_{\vec{w}})})\\
	&\leq& C\| \vec{b} \|_{BMO^{m}}\|\mathcal{M}_{s}(\vec{f})\|_{L^{p}(v_{\vec{w}})}=C\| \vec{b} \|_{BMO^{m}}\|\mathcal{M}(\vec{f^{s}})\|_{L^{\frac{p}{s}}(v_{\vec{w}})}^{\frac{1}{s}} \\
	&\leq & C\| \vec{b} \|_{BMO^{m}}\prod_{j=1}^{m}\||f_{j}|^{s}\|_{L^{\frac{p_{j}}{s}(w_{j})}}^{\frac{1}{s}}=C\| \vec{b} \|_{BMO^{m}}\prod_{j=1}^{m}\|f_{j}\|_{L^{p_{j}}(w_{j})}.
\end{eqnarray*}
This finishes the proof of the Theorem 3.2.$\hfill\square$

\medskip
\noindent{{\it Proof of Theorem $3.3$.}}
For the sake of simplicity, we only consider the case $m=2$. The proof of other cases is similar.
Let $f_{1}$, $f_{2}$ be bounded measurable functions with compact support and $b_{1}$, $b_{2}\in BMO$. Then for any constants $\lambda_{1}$ and $\lambda_{2},$
\begin{eqnarray*}
T_{\prod\vec{b}}(\vec{f})(z)&=&(b_{1}(z)-\lambda_{1})(b_{2}(z)-\lambda_{2})T(f_{1},f_{2})(z)-(b_{1}(z)-\lambda_{1}))T(f_{1},(b_{2}-\lambda_{2})f_{2})(z)\\
&&-(b_{2}(z)-\lambda_{2})T((b_{1}-\lambda_{1})f_{1},f_{2})(z)+T((b_{1}-\lambda_{1})f_{1},(b_{2}-\lambda_{2})f_{2})(z)\\
&=&-(b_{1}(z)-\lambda_{1})(b_{2}(z)-\lambda_{2})T(f_{1},f_{2})(z)+(b_{1}(z)-\lambda_{1})T^{2}_{b_{2}-\lambda_{2}}(f_{1},f_{2})(z)\\
&&+(b_{2}(z)-\lambda_{2})T^{1}_{b_{1}-\lambda_{1}}(f_{1},f_{2})(z)+T((b_{1}-\lambda_{1})f_{1},(b_{2}-\lambda_{2})f_{2})(z).
\end{eqnarray*}

Let $C_{0}$ be a constant determined later. Fixed $x\in \mathbb{R}^{n}$, for any ball $B(x_{0},r_{B})$
containing $x$ and $0<\delta<\frac{1}{2}$, we have
\begin{eqnarray*}
	&&\biggl(\frac{1}{|B|}\dint_{B}||T_{\prod\vec{b}}(\vec{f})(z)|^{\delta}-|C_{0}|^{\delta}|dz
	\biggr)^{\frac{1}{\delta}}\\
	&\leq&\biggl(\frac{1}{|B|}\dint_{B}|T_{\prod\vec{b}}(\vec{f})(z)-C_{0}|^{\delta}dz
	\biggr)^{\frac{1}{\delta}}\\
	&\leq&C\biggl(\frac{1}{|B|}\dint_{B}|(b_{1}(z)-\lambda_{1})(b_{2}(z)-\lambda_{2})T(f_{1},f_{2})(z)|^{\delta}dz
	\biggr)^{\frac{1}{\delta}}\\
	&&+C\biggl(\frac{1}{|B|}\dint_{B}|(b_{1}(z)-\lambda_{1})T^{2}_{b_{2}-\lambda_{2}}(f_{1},f_{2})(z)|^{\delta}dz
	\biggr)^{\frac{1}{\delta}}\\
	&&+C
	\biggl(\frac{1}{|B|}\dint_{B}|(b_{2}(z)-\lambda_{2})T^{1}_{b_{1}-\lambda_{1}}(f_{1},f_{2})(z)|^{\delta}dz
	\biggr)^{\frac{1}{\delta}}\\
	&&+C
	\biggl(\frac{1}{|B|}\dint_{B}|T((b_{1}-\lambda_{1})f_{1},(b_{2}-\lambda_{2})f_{2})(z)-C_{0}|^{\delta}dz
	\biggr)^{\frac{1}{\delta}}\\
	&:=&I+II+III+IV.
\end{eqnarray*}

Then we calculate each term respectively. 
Denote $B^{*}=16\sqrt{2}B$ and $\lambda_{j}=(b_{j})_{B^{*}}=\frac{1}{|16\sqrt{2}B|}\int_{16\sqrt{2}B}b_{j}(x)dx$, $j=1,2$. Since $0<\delta<\frac{1}{2}$ and $0<\delta<\varepsilon_{0}<\infty$,
it is easy to seek an $l$ such that $1<l<\min\{\frac{\varepsilon_{0}}{\delta}$, $\frac{1}{1-\delta}\}$. We can deduce that $l\delta<\varepsilon_{0}$ and $l'\delta>1$. Choose $q_{1},q_{2}\in (1,\infty)$ satisfying  $\frac{1}{q_{1}}+\frac{1}{q_{2}}=\frac{1}{l'}$, then $\frac{1}{q_{1}}+\frac{1}{q_{2}}+\frac{1}{l}=1$,  $q_{1}\delta>1$ and $q_{2}\delta>1$. By Lemma 2.4, we have
\begin{eqnarray*}
	I&\leq&C\biggl(\frac{1}{|B|}\dint_{B}|T(f_{1},f_{2})(z)|^{l\delta }dz\biggr)^{\frac{1}{l\delta}}\biggl(\frac{1}{|B|}\dint_{B}|b_{1}(z)-\lambda_{1}|^{q_{1}\delta}dz\biggr)^{\frac{1}{q_{1}\delta}}
    \\
	&&\times \biggl(\frac{1}{|B|}\dint_{B}|b_{2}(z)-\lambda_{2}|^{q_{2}\delta }dz\biggr)^{\frac{1}{q_{2}\delta }}\\
	&\leq&C\|b_{1}\|_{BMO}\|b_{2}\|_{BMO}M_{\delta l}(T(f_{1},f_{2}))(x)\\
	&\leq&C\|b_{1}\|_{BMO}\|b_{2}\|_{BMO}M_{\varepsilon_{0}}(T(f_{1},f_{2}))(x).
	\end{eqnarray*}

It follow from H\"{o}lder's inequality that
\begin{eqnarray*}
	II&\leq&C\biggl(\frac{1}{|B|}\dint_{B}|T_{b_{2}-\lambda_{2}}^{2}(f_{1},f_{2})(z)|^{l\delta}dz\biggr)^{\frac{1}{l\delta}}\biggl(\frac{1}{|B|}\dint_{B}|b_{1}(z)-\lambda_{1}|^{l'\delta}dz\biggr)^{\frac{1}{l'\delta}}\\
	&\leq&C\|b_{1}\|_{BMO}M_{\delta l}(T_{b_{2}-\lambda_{2}}^{2}(f_{1},f_{2}))(x)\\
	&\leq& C\|b_{1}\|_{BMO}M_{\varepsilon_{0}}(T_{b_{2}}^{2}(f_{1},f_{2}))(x).
	\end{eqnarray*}
Similarly we can get
$$III\leq C\|b_{2}\|_{BMO}M_{\varepsilon_{0}}(T_{b_{1}}^{1}(f_{1},f_{2}))(x).$$

In the next, we discuss $IV$. Split $f_{i}$ into two parts
$f_{i}=f_{i}^{0}+f_{i}^{\infty}$,
where $f_{i}^{0}=f\chi_{B^{*}}$
and $f_{i}^{\infty}=f_{i}-f_{i}^{0}$, $i=1, 2$. Choose $z_{0}\in4B\backslash3B$.  Denote $\vec{z_{0}}=(z_{0},z_{0})$ and
$\vec{y}=(y_{1},y_{2})$. Let
$$C_{0}=\sum_{(\alpha_{1},\alpha_{2})\in\Gamma}T((b_{1}-\lambda_{1})f_{1}^{\alpha_{1}},(b_{2}-\lambda_{2})f_{2}^{\alpha_{2}})(z_{0}),$$
where $\Gamma=\{(\alpha_{1},\alpha_{2})$: there is at least one $\alpha_{j}\neq 0,j=1,2\}$, then
\begin{eqnarray*}
	IV&\leq&C\biggl(\frac{1}{|B|}\dint_{B}|T((b_{1}-\lambda_{1})f_{1}^{0},(b_{2}-\lambda_{2})f_{2}^{0})(z)|^{\delta}dz
	\biggr)^{\frac{1}{\delta}}\\
	&&+C\sum_{(\alpha_{1},\alpha_{2})\in\Gamma}\biggl(\frac{1}{|B|}\dint_{B}|T((b_{1}-\lambda_{1})f_{1}^{\alpha_{1}},(b_{2}-\lambda_{2})f_{2}^{\alpha_{2}})(z)\\
	&&-T((b_{1}-\lambda_{1})f_{1}^{\alpha_{1}},(b_{2}-\lambda_{2})f_{2}^{\alpha_{2}})(z_{0})|^{\delta}dz
	\biggr)^{\frac{1}{\delta}}\\
	&:=&IV_{0}+\sum_{(\alpha_{1},\alpha_{2})\in\Gamma}IV_{\alpha_{1},\alpha_{2}}.
\end{eqnarray*}

 Let $t=s/q'$, then it follows from $s>q'$ that $t>1$. Due to for any $j=1,2$, $r_{j}\leq q'<s$, we can get $r_{j}t\leq s$.
 By Lemma 2.3 and H\"{o}lder's inequality,
\begin{eqnarray*}
	IV_{0}&\leq& C|B|^{-1/\delta}\|T((b_{1}-\lambda_{1})f_{1}^{0},(b_{2}-\lambda_{2})f_{2}^{0})\|_{L^{\delta}(B)} \\ &\leq&C|B|^{-1/r}\|T((b_{1}-\lambda_{1})f_{1}^{0},(b_{2}-\lambda_{2})f_{2}^{0})\|_{L^{r,\infty}(B)} \\
	&\leq&C \biggl(\frac{1}{|16\sqrt{2}B|}\dint_{16\sqrt{2}B}|b_{1}(y_{1})-\lambda_{1}|^{r_{1}}|f_{1}(y_{1})|^{r_{1}}dy_{1}\biggr)^{\frac{1}{r_{1}}}\\
	&&\times\biggl(\frac{1}{|16\sqrt{2}B|}\dint_{16\sqrt{2}B}|b_{2}(y_{2})-\lambda_{2}|^{r_{2}}|f_{2}(y_{2})|^{r_{2}}dy_{2}\biggr)^{\frac{1}{r_{2}}}\\
	&\leq&C \biggl(\frac{1}{|16\sqrt{2}B|}\dint_{16\sqrt{2}B}|f_{1}(y_{1})|^{r_{1}t}dy_{1}\biggr)^{\frac{1}{r_{1}t}}\biggl(\frac{1}{|16\sqrt{2}B|}\dint_{16\sqrt{2}B}|b_{1}(y_{1})-\lambda_{1}|^{r_{1}t'}dy_{1}\biggr)^{\frac{1}{r_{1}t'}}\\
	&&\times \biggl(\frac{1}{|16\sqrt{2}B|}\dint_{16\sqrt{2}B}|f_{2}(y_{2})|^{r_{2}t}dy_{2}\biggr)^{\frac{1}{r_{2}t}}\biggl(\frac{1}{|16\sqrt{2}B|}\dint_{16\sqrt{2}B}|b_{2}(y_{2})-\lambda_{2}|^{r_{2}t'}dy_{2}\biggr)^{\frac{1}{r_{2}t'}}\\
	&\leq&C\|b_{1}\|_{BMO}\|b_{2}\|_{BMO}\biggl(\frac{1}{|16\sqrt{2}B|}\dint_{16\sqrt{2}B}|f_{1}(y_{1})|^{s}dy_{1}\biggr)^{\frac{1}{s}}\\
	&&\times\biggl(\frac{1}{|16\sqrt{2}B|}\dint_{16\sqrt{2}B}|f_{2}(y_{2})|^{s}dy_{2}\biggr)^{\frac{1}{s}}\\
	&\leq&C\|b_{1}\|_{BMO}\|b_{2}\|_{BMO}\mathcal{M}_{s}(\vec{f})(x).
\end{eqnarray*}

For any $(\alpha_{1},\alpha_{2}) \in \Gamma$, we can find a $ j\in\{1, 2\}$ that makes  $\alpha_{j}=\infty$. Then for any $\vec{y}=(y_{1},y_{2})\in {\rm supp}f_{1}^{\alpha_{1}}\times {\rm supp}f_{2}^{\alpha_{2}}$ and $z\in B$, 
$|\vec{y}-\vec{z_{0}}|\geq 2\sqrt{2}|z-z_{0}|$ and
$2r_{B}\leq|z-z_{0}|\leq5r_{B}$.
Thus
\begin{eqnarray*}
	IV_{\alpha_{1},\alpha_{2}}
	&\leq& C\frac{1}{|B|}\dint_{B}\dint_{|\vec{y}-\vec{z_{0}}|\geq 2\sqrt{2}|z-z_{0}|}|K(z,y_{1},y_{2})-K(z_{0},y_{1},y_{2})|\\
	&&\times\prod_{j=1}^{2}|b_{j}(y_{j})-\lambda_{j}||f_{j}(y_{j})|d\vec{y}dz  
	\\
	&\leq&C\frac{1}{|B|}\dint_{B}\sum_{k=1}^{\infty}\dint_{2^{k}\sqrt{2}|z-z_{0}|\leq|\vec{y}-\vec{z_{0}}|< 2^{k+1}\sqrt{2}|z-z_{0}|}|K(z,y_{1},y_{2})-K(z_{0},y_{1},y_{2})|\\
	&&\times\prod_{j=1}^{2}|b_{j}(y_{j})-\lambda_{j}||f_{j}(y_{j})|d\vec{y}dz 
	\\
	&\leq&C\frac{1}{|B|}\dint_{B}\sum_{k=1}^{\infty}\biggl(\dint_{2^{k}\sqrt{2}|z-z_{0}|\leq|\vec{y}-\vec{z_{0}}|< 2^{k+1}\sqrt{2}|z-z_{0}|}|K(z,y_{1},y_{2})-K(z_{0},y_{1},y_{2})|^{q}d\vec{y}\biggr)^{\frac{1}{q}}\\
	&&\times\biggl(\dint_{|\vec{y}-\vec{z_{0}}|<2^{k+1}\sqrt{2}|z-z_{0}|}\prod_{j=1}^{2}|b_{j}(y_{j})-\lambda_{j}|^{q'}|f_{j}(y_{j})|^{q'}
	d\vec{y}\biggr)^{\frac{1}{q'}}dz 
	\\
	&\leq&C\frac{1}{|B|}\dint_{B}\sum_{k=1}^{\infty}\biggl(\dint_{2^{k}\sqrt{2}|z-z_{0}|\leq|\vec{y}-\vec{z_{0}}|< 2^{k+1}\sqrt{2}|z-z_{0}|}|K(z,y_{1},y_{2})-K(z_{0},y_{1},y_{2})|^{q}d\vec{y}\biggr)^{\frac{1}{q}}\\
	&&\times\prod_{j=1}^{2}\biggl(\frac{1}{|B(z_{0},2^{k+1}\sqrt{2}|z-z_{0}|)|}\dint_{|y_{j}-z_{0}|< 2^{k+1}\sqrt{2}|z-z_{0}|}|b_{j}(y_{j})-\lambda_{j}|^{q'}|f_{j}(y_{j})|^{q'}dy_{j}\biggr)^{\frac{1}{q'}}\\
	&&\times  |B(z_{0},2^{k+1}\sqrt{2}|z-z_{0}|)|^{\frac{2}{q'}}dz\\
	&\leq&C\frac{1}{|B|}\dint_{B}\sum_{k=1}^{\infty}\biggl(\dint_{2^{k}\sqrt{2}|z-z_{0}|\leq|\vec{y}-\vec{z_{0}}|< 2^{k+1}\sqrt{2}|z-z_{0}|}|K(z,y_{1},y_{2})-K(z_{0},y_{1},y_{2})|^{q}d\vec{y}\biggr)^{\frac{1}{q}}\\
	&&\times\prod_{j=1}^{2}\biggl(\frac{1}{|B(x_{0},2^{k+2}\sqrt{2}|z-z_{0}|)|}\dint_{|y_{j}-x_{0}|< 2^{k+2}\sqrt{2}|z-z_{0}|}|b_{j}(y_{j})-\lambda_{j}|^{q't'}|dy_{j}\biggr)^{\frac{1}{q't'}}\\
	&&\times\prod_{j=1}^{2}\biggl(\frac{1}{|B(z_{0},2^{k+1}\sqrt{2}|z-z_{0}|)|}\dint_{|y_{j}-z_{0}|< 2^{k+1}\sqrt{2}|z-z_{0}|}|f_{j}(y_{j})|^{q't}dy_{j}\biggr)^{\frac{1}{q't}} \\
	&&\times|B(z_{0},2^{k+1}\sqrt{2}|z-z_{0}|)|^{\frac{2}{q'}}dz\\
	&\leq&C\|b_{1}\|_{BMO}\|b_{2}\|_{BMO}\frac{1}{|B|}\dint_{B}
	\sum_{k=1}^{\infty}|z-z_{0}|^{- \frac{2n}{q'}}C_{k}
	2^{-\frac{2kn}{q'}}k^{2}(2^{k+1}\sqrt{2}|z-z_{0}|)^{ \frac{2n}{q'}}\\
	&&\times\prod_{j=1}^{2}\biggl(\frac{1}{|B(z_{0},2^{k+1}\sqrt{2}|z-z_{0}|)|}\dint_{|y_{j}-z_{0}|< 2^{k+1}\sqrt{2}|z-z_{0}|}
	|f_{j}(y_{j})|^{s}dy_{j}\biggr)^{\frac{1}{s}}dz\\
	&\leq& C\|b_{1}\|_{BMO}\|b_{2}\|_{BMO}\mathcal{M}_{s}(\vec{f})(x).
\end{eqnarray*}
Thus
$$
	\sum_{(\alpha_{1},\alpha_{2})\in\Gamma}IV_{\alpha_{1},\alpha_{2}}
	\leq C\|b_{1}\|_{BMO}\|b_{2}\|_{BMO}\mathcal{M}_{s}(\vec{f})(x).
$$
	So
	$$ IV\leq C\|b_{1}\|_{BMO}\|b_{2}\|_{BMO}\mathcal{M}_{s}(\vec{f})(x).$$
	
	Finally, we can get the result,
\begin{eqnarray*}
	M^{\sharp}_{\delta}(T_{\prod\vec{b}}(\vec{f}))(x)&=&	M^{\sharp}(|T_{\prod\vec{b}}(\vec{f})|^{\delta})^{\frac{1}{\delta}}(x)\\
	&\leq&C\dsup_{B\ni x}\biggl(\frac{1}{|B|}\dint_{B}||T_{\prod\vec{b}}(\vec{f})(z)|^{\delta}-|C_{0}|^{\delta}|dz
	\biggr)^{\frac{1}{\delta}}\\
	&\leq& C \|b_{1}\|_{BMO}\|b_{2}\|_{BMO}\biggl(M_{\varepsilon_{0}}(T(\vec{f}))(x)+\mathcal{M}_{s}(\vec{f})(x)\biggr)\\
	&&+C(\|b_{1}\|_{BMO}M_{\varepsilon_{0}}(T_{b_{2}}^{2}(\vec{f}))(x)
	+\|b_{2}\|_{BMO}M_{\varepsilon_{0}}(T_{b_{1}}^{1}(\vec{f}))(x)).
\end{eqnarray*}
This completes the proof of Theorem 3.3.$\hfill\square$

\medskip
\noindent{{\it Proof of Theorem $3.4$.}}

Proceeding as in the proof of Theorem 3.2, we can find 
a $s$ such that $q'<s$, $s<p_{j}$, $j=1,\dots,m$ and
$\vec{w}\in A_{\vec{P}/s}$.

Choose $\delta,\varepsilon_{1},\varepsilon_{2},\dots,\varepsilon_{m}$
satisfing $0<\delta<\varepsilon_{1}<\varepsilon_{2}<\dots<\varepsilon_{m}<\frac{1}{m}$. We can deduce the conclusion by Lemma 2.5 and Lemma 2.8.
$$\|M_{\varepsilon_{j}}(T(\vec{f}))\|_{L^{p}(v_{\vec{w}})}
\leq C\|M_{\varepsilon_{j}}^{\sharp}(T(\vec{f}))\|_{L^{p}(v_{\vec{w}})}
\leq C\|\mathcal{M}_{q'}(\vec{f})\|_{L^{p}(v_{\vec{w}})} 
\leq C\|\mathcal{M}_{s}(\vec{f})\|_{L^{p}(v_{\vec{w}})}.$$
By Theorem 3.3, we have
\begin{eqnarray*}
	\|M^{\sharp}_{\delta}(T_{\prod\vec{b}}(\vec{f}))\|_{L^{p}(v_{\vec{w}})}
	&\leq& C\prod_{j=1}^{m}\| b_{j} \|_{BMO} \biggl(\|M_{\varepsilon_{1}}(T(\vec{f}))\|_{L^{p}(v_{\vec{w}})}+\|\mathcal{M}_{s}(\vec{f})\|_{L^{p}(v_{\vec{w}})}\biggr)\\
	&&+C\sum_{j=1}^{m-1}\sum_{\sigma\in C_{j}^{m}}\prod_{i=1}^{j}\| b_{\sigma(i)} \|_{BMO} \|M_{\varepsilon_{1}}(T_{\prod\vec{b}_{\sigma' }}(\vec{f}))\|_{L^{p}(v_{\vec{w}})}\\
	&\leq& C\prod_{j=1}^{m}\| b_{j} \|_{BMO} \biggl(\|M_{\varepsilon_{1}}(T(\vec{f}))\|_{L^{p}(v_{\vec{w}})}+\|\mathcal{M}_{s}(\vec{f})\|_{L^{p}(v_{\vec{w}})}\biggr)\\
	&&+C\sum_{j=1}^{m-1}\sum_{\sigma\in C_{j}^{m}}\prod_{i=1}^{j}\| b_{\sigma(i)} \|_{BMO} \|M_{\varepsilon_{1}}^{\sharp}(T_{\prod\vec{b}_{\sigma' }}(\vec{f}))\|_{L^{p}(v_{\vec{w}})}.
\end{eqnarray*}

For the purpose of reducing the dimension of BMO functions in the commutators, we apply the Theorem 3.3 to $\|M_{\varepsilon_{1}}^{\sharp}(T_{\prod\vec{b}_{\sigma' }}(\vec{f}))(x))\|_{L^{p}(v_{\vec{w}})}.$

Let $\sigma=\{\sigma(1),\dots,\sigma(j)\}$ and $\sigma'=\{\sigma(j+1),\dots,\sigma(m)\}$, $A_{h}$=\{$\sigma_{1}$: any finite subset of $\sigma'$ with different elements\} and $\sigma_{1}'=\sigma'-\sigma_{1}$,

\begin{eqnarray*}
	\|M^{\sharp}_{\varepsilon_{1}}(T_{\prod\vec{b}_{\sigma'}}(\vec{f}))\|_{L^{p}(v_{\vec{w}})}
	&\leq& C\prod_{i=j+1}^{m}\| b_{\sigma(i)} \|_{BMO} \biggl(\|M_{\varepsilon_{2}}(T(\vec{f}))\|_{L^{p}(v_{\vec{w}})}+\|\mathcal{M}_{s}(\vec{f})\|_{L^{p}(v_{\vec{w}})}\biggr)\\
	&&+C\sum_{h=1}^{m-j-1}\sum_{\sigma_{1}\in A_{h}}\prod_{i=1}^{h}\| b_{\sigma_{1}(i)}\|_{BMO} \|M_{\varepsilon_{2}}(T_{\prod\vec{b}_{\sigma_{1}' }}(\vec{f}))\|_{L^{p}(v_{\vec{w}})}.\\
\end{eqnarray*}
Repeating the process above and using Theorem 3.1, we can get
\begin{eqnarray*}
	\|M^{\sharp}_{\delta}(T_{\prod\vec{b}}(\vec{f}))\|_{L^{p}(v_{\vec{w}})}
	&\leq& C\prod_{j=1}^{m}\| b_{j} \|_{BMO} \biggl(A_{m+1}(m,n)\|\mathcal{M}_{s}(\vec{f})\|_{L^{p}(v_{\vec{w}})}+A_{1}(m,n)\\
	&&\times\|M_{\varepsilon_{1}}(T(\vec{f}))\|_{L^{p}(v_{\vec{w}})}
	+A_{2}(m,n)\|M_{\varepsilon_{2}}(T(\vec{f}))\|_{L^{p}(v_{\vec{w}})}+\cdots \\
	&&+A_{m}(m,n)\|M_{\varepsilon_{m}}(T(\vec{f}))\|_{L^{p}(v_{\vec{w}})}\biggr),\\
	\end{eqnarray*}
where $A_{1}(m,n),A_{2}(m,n),\dots,A_{m+1}(m,n)$ are finite real numbers related to $m$ and $n$.

We can deduce the following conclusion from Lemma 2.5, Lemma 2.7 and $\vec{w}\in A_{\vec{P}/s}$.
\begin{eqnarray*}
	\|T_{\prod\vec{b}}(\vec{f})\|_{L^{p}(v_{\vec{w}})}
	&\leq&\|M_{\delta}(T_{\prod\vec{b}}(\vec{f}))\|_{L^{p}(v_{\vec{w}})}\\
	&\leq&C\|M^{\sharp}_{\delta}(T_{\prod\vec{b}}(\vec{f}))\|_{L^{p}(v_{\vec{w}})}\\
	&\leq& C\prod_{j=1}^{m}\| b_{j} \|_{BMO} \biggl(A_{m+1}(m,n)\|\mathcal{M}_{s}(\vec{f})\|_{L^{p}(v_{\vec{w}})}+A_{1}(m,n)\\
	&&\times\|M_{\varepsilon_{1}}(T(\vec{f}))\|_{L^{p}(v_{\vec{w}})}+A_{2}(m,n)\|M_{\varepsilon_{2}}(T(\vec{f}))\|_{L^{p}(v_{\vec{w}})}+\dots \\
	&&+A_{m}(m,n)\|M_{\varepsilon_{m}}(T(\vec{f}))\|_{L^{p}(v_{\vec{w}})}\biggr)\\	
	&\leq& C \prod_{j=1}^{m}\| b_{j} \|_{BMO}\|\mathcal{M}_{s}(\vec{f})\|_{L^{p}(v_{\vec{w}})}\\
	&\leq& C \prod_{j=1}^{m}\| b_{j} \|_{BMO}\|\|f_{j}\|_{L^{p_{j}}(w_{j})}.
\end{eqnarray*}
The proof is completed.$\hfill\square$

\medskip
\noindent{{\bf Remark 3.5.}}  It should be pointed out that the corresponding results of multilinear commutators and iterated commutators of multilinear singular integral operators with Dini's type kernels in \cite{lz1,sz1,sz2} and the standard kernels can be deduced as special cases of our results in this paper.

{\centering
\medskip

LIWEN GAO\par
SCHOOL OF SCIENCE\par
CHINA UNIVERSITY OF MINING AND TECHNOLOGY, BEIJING\par
BEIJING 100083, P. R. CHINA \par
Email address: {\bf gaoliwen1206@163.com}\par
\centerline{}
YAN LIN\par
SCHOOL OF SCIENCE\par
CHINA UNIVERSITY OF MINING AND TECHNOLOGY, BEIJING\par
BEIJING 100083, P. R. CHINA \par
Email address: {\bf linyan@cumtb.edu.cn}\par
\centerline{}
SHUHUI YANG\par
SCHOOL OF SCIENCE\par
CHINA UNIVERSITY OF MINING AND TECHNOLOGY, BEIJING\par
BEIJING 100083, P. R. CHINA \par
Email address: {\bf yangshuhui0601@163.com}

	\end{document}